\documentstyle[12pt,amscd]{amsart}
\newtheorem{pr}{Proposition}
\newtheorem{lm}{Lemma}
\newtheorem{tm}{Theorem}
\newtheorem{cor}{Corollary}
\newcommand{\proj}{\bold P}
\newcommand{\grass}{\bold G}
\newcommand{\barr}{\overline}

\newcommand{\rarr}{\rightarrow}
\newcommand{\oh}{{\cal{O}}}

\newcommand{\eqq}{\stackrel{\sim}{=}}

\newcommand{\com}{\Bbb{C}}
\newcommand{\Q}{{\Bbb Q}}

\begin{document}
\title{The Chow Ring of the Non-Linear Grassmannian}
\author{ Rahul Pandharipande$^1$}
\date{11 March 1996}
\maketitle
\pagestyle{plain}
\footnotetext[1]{Research 
partially supported by an NSF Post-Doctoral Fellowship.}
\setcounter{section}{-1}
\section{$\bold{Summary}$}
\label{intro}
Let $\com$ be the ground field of complex numbers.
Let $1\leq k \leq r$ be integers. 
The Grassmannian  $\grass(\proj^k,\proj^r)$ of
projective $k$-planes in $\proj^r$ can be viewed as
the moduli space of (unparameterized) regular maps from
$\proj^k$ to $\proj^r$ of degree $1$. 
Let $M_{\proj^k}(\proj^r,d)$ be the 
coarse moduli space of (unparameterized) regular
maps $\mu:\proj^k \rarr \proj^r$
satisfying $\mu^{*}(\oh_{\proj^r}(1))\eqq \oh_{\proj^k}(d)$.
Two maps 
$$\mu:\proj^k \rarr \proj^r, \ \mu':\proj^k \rarr \proj^r$$
are equivalent for the moduli problem if
there exists an element
$\sigma\in \bold{PGL}_{k+1}$ satisfying
$\mu' \circ \sigma = \mu$. 
If $\mu:\proj^k \rarr \proj^r$ is a non-constant 
regular map, it is easy to
show that $dim(Im(\mu))=k$ and $\mu: \proj^k \rarr Im(\mu)$
is a {\em finite} morphism. The space $M_{\proj^k}(\proj^r,d)$
is a natural non-linear generalization of the Grassmannian.

In section (\ref{construct}), 
$M_{\proj^k}(\proj^r,d)$ will be constructed
via Geometric Invariant Theory. $M_{\proj^k}(\proj^r,d)$
is an irreducible, normal, quasi-projective variety
with finite quotient singularities. 
Let $A_i\big(M_{\proj^k}(\proj^r,d)\big)
\otimes {\Bbb Q}$
be the Chow group (tensor ${\Bbb Q}$)
 of $i$-cycles modulo linear equivalence.
Since the space $M_{\proj^k}(\proj^r,d) $
has finite quotient singularities, the Chow groups $\bigoplus
 (A_i\otimes
{\Bbb Q})$
naturally form a graded ring via intersection. Since
$\Q$-coefficients are required for the intersection theory,
all Chow groups considered here will be taken with
$\Q$-coefficients. 
Let $Ch(k,r,d)$ denote the Chow ring of $M_{\proj^k}(\proj^r,d) $.
The ring $Ch(k,r,1)$ is simply the Chow ring of the
linear Grassmannian $\grass(\proj^k,\proj^r)$. 
The main result of this paper is a determination of
$Ch(k,r,d)$.
\begin{tm}
\label{mainn}
There is a {\em canonical} isomorphism of graded rings 
$$\lambda: Ch(k,r,d) \rarr Ch(k,r,1).$$
\end{tm}

Let $\barr{M}_{0,n}(\proj^r,d)$ be the coarse
moduli space of $n$-pointed Kontsevich stable maps from
a genus $0$ curve to $\proj^r$. 
Let ${M}_{0,n}(\proj^r,d) \subset \barr{M}_{0,n}(\proj^r,d) $
denote the non-empty open set corresponding to 
$n$-pointed, stable maps from $\proj^1$ to $\proj^r$.
The complement of ${M}_{0,n}(\proj^r,d)$ in
$\barr{M}_{0,n}(\proj^r,d)$ consists of the stable
maps with reducible domains. 
A foundational treatment of these moduli spaces of pointed
stable maps of genus $0$ curves can be found in [K], [KM], and [FP].
The spaces $M_{0,0}(\proj^r,d)$ and $M_{\proj^1}(\proj^r,d)$
are identical. The following corollary
is therefore a special case of Theorem (\ref{mainn}).
\begin{cor}
\label{corr}
The Chow ring (with $\Q$-coefficients) of 
$M_{0,0}(\proj^r,d)$ is canonically isomorphic to
the Chow ring of the Grassmannian $\grass(\proj^1, \proj^r)$.
\end{cor}

Corollary
 (\ref{corr}) is related by loose analogy to
results and conjectures on the Chow ring of $M_g$.
C. Faber has studied the subring of the Chow ring
of $M_g$ generated by certain geometric classes.
Faber has conjectured a presentation of this subring (which
may be the entire Chow ring of $M_g$). The conjectured
ring looks like the cohomology ring of a compact
manifold -- for example, it satisfies Poincar\'e duality.
$M_{0,0}(\proj^r,d)\subset 
\barr{M}_{0,0}(\proj^r,d)$ is a  zero-pointed open cell
analogous to $M_g\subset \barr{M}_g$. Corollary (\ref{corr}), then, is
analogous to Faber's conjectures.

In [GP], the Poincar\'e polynomial of 
$\barr{M}_{0,n}(\proj^r,d)$ is computed.
The {\em virtual} Poincar\'e polynomial of
${M}_{0,0}(\proj^r,d)$ is needed as a preliminary
result. It was found the virtual Poincar\'e polynomial
of ${M}_{0,0}(\proj^r,d)$ is essentially the Poincar\'e
polynomial of $\bold{G}(\proj^1, \proj^r)$. This observation
provided the starting point for Theorem (\ref{mainn}).
Thanks are especially due to E. Getzler for many discussions
about the geometry of the space $M_{0,0}(\proj^r,d)$.
The theory of
equivariant Chow groups ([EG], [T]) plays an essential role
in the proof of Theorem (\ref{mainn}). The author wishes
to thank D. Edidin, W. Graham, and B. Totaro for
the long discussions in which this theory was explained.
The author has also benefitted from conversations with
P. Belorouski, W. Fulton, and H. Tamvakis.

\section{$M_{\proj^k}(\proj^r,d)$}
\label{construct}
A family of degree $d$ maps of $\proj^k$ to $\proj^r$ consists of the
data $(\pi:\cal{P}\rarr S, \ \mu: \cal{P}\rarr \proj^r)$
where:
\begin{enumerate}
\item[(i.)] 
$S$ is a noetherian scheme of finite type over $\com$.
\item[(ii.)]  $\pi:\cal{P}\rarr S$ is a flat projective morphism
with geometric fibers isomorphic to $\proj^k$.
\item[(iii.)] The restriction of $\mu^{*}(\oh_{\proj^r}(1))$
to each geometric fiber of $\pi$ is isomorphic to
 $\oh_{\proj^k}(d)$.
\end{enumerate} 
Two families of maps over $S$,  
$$(\pi:\cal{P}\rarr S, \ \mu), \ (\pi':\cal{P}'\rarr S, \ \mu' )$$ are
isomorphic if there exists an isomorphism of 
schemes $\sigma: \cal{P} \rarr \cal{P}'$
such that
$$\mu=\mu' \circ \sigma, \ \pi = \pi' \circ \sigma.$$ 
Let $\cal{M}_{\proj^k}(\proj^r,d)$ be the
contravariant functor from schemes to sets defined
as follows. $\cal{M}_{\proj^k}(\proj^r,d) \ (S)$ is the
set of isomorphism classes  of families over $S$ of
degree $d$ maps from $\proj^k$  to $\proj^r$.
 
A coarse moduli space $M_{\proj^k}(\proj^r,d)$
is easily obtained via Geometric Invariant Theory.
Care is taken
here to exhibit $M_{\proj^k}(\proj^r,d)$ as a quotient
of a proper $\bold{GL}_{k+1}$-action with finite stabilizers. 
In section
(\ref{genarg}), 
the equivariant Chow groups of this $\bold{GL}_{k+1}$-action
will be analyzed.

Let
$$U(k,r,d) \subset \bigoplus_{0}^{r}H^0(\proj^k, \oh_{\proj^k}
(d))$$
be the Zariski open locus of basepoint free
$(r+1)$-tuples of polynomials. There is
a natural $\bold{GL}_{k+1}$-action on
$\bigoplus_{0}^{r}H^0(\proj^k, \oh_{\proj^k}
(d))$ obtained from the naturally linearized action
of $\bold{GL}_{k+1}$ on $\proj^k$. 
This $\bold{GL}_{k+1}$-action
leaves
$U(k,r,d)$ invariant. Note, since every
regular map $\mu: \proj^k \rarr \proj^r$ is
finite onto its image,  $\bold{GL}_{k+1}$
acts with {\em finite stabilizers} on $U(k,r,d)$.

Let $\bold{1}\eqq \com$ be a $1$ dimensional complex vector
space with the trivial $GL_{k+1}$-action.
Let $Det$ be the $1$ dimensional determinant representation
of $GL_{k+1}$. For convenience, let $Z$ denote
$\bigoplus_{0}^r
H^0(\proj^k, \oh_{\proj^k}(d))$.
There is a $GL_{k+1}$-equivariant  inclusion
$$U(k,r,d) \subset \proj\big( Det \otimes\ 
(\bold{1} \oplus   Sym^q( Z) \oplus Z ) \big)$$
obtained by the following equation:
\begin{equation}
\label{contort}
\xi \in U(k,r,d) \rarr [\ 1\otimes\  
(1\oplus (\xi\otimes \cdots \otimes \xi) \oplus \xi)\ ].
\end{equation}
The representations $Det$ and $Sym^q(Z)$  in (\ref{contort})
occur to obtain the correct
G.I.T. linearization. The final $Z$ factor occurs
to insure (\ref{contort}) is an inclusion (consider
scaling  $U(k,r,d)$ by a constant
$q^{th}$-root of unity).  
\begin{lm}
\label{geintth}
Consider the naturally linearized 
action of $GL_{k+1}$ on $$ 
\proj\big( Det \otimes\ 
(\bold{1} \oplus   Sym^q( Z) \oplus Z ) \big).$$
Then, for $q>k+1$,
$$U(k,r,d)  \subset 
\proj\big( Det \otimes\ 
(\bold{1} \oplus   Sym^q( Z) \oplus Z ) \big)
^{stable}.$$
\end{lm}
\begin{pf}
The Lemma is a consequence of the
Numerical Criterion of stability.
A development of Geometric Invariant Theory
can be found in [MFK] and [N]. Let
$V_{k+1}$ be a $(k+1)$-dimensional $\com$-vector
space such that $\proj^k= \proj (V_{k+1})$.
Let $\barr{v}=v_0, \ldots, v_{k}$ be a basis 
of $V_{k+1}$ with 
integer weights $w_0, \ldots, w_{k}$ (not all zero).
Let $\xi \in  U(k,r,d)$ correspond to a basepoint
free map determined (in the basis $\barr{v}$)
by $[f_0, \ldots , f_r]$ where
each $f_l$ is an an element of $Sym^d(V_{k+1}^*)$.
The diagonal coordinates of 
$$\xi \in 
\proj\big( Det \otimes\ 
(\bold{1} \oplus   Sym^q( Z) \oplus Z ) \big)$$
with respect to the $\com^*$-action determined by the
weights and basis $\barr{v}$ are the following:
$$ 1\otimes 1 \in Det\otimes \bold{1},$$
\begin{equation}
\label{termmm}
 1 \otimes  (\xi \otimes \cdots \otimes 
\xi) \in Det \otimes Sym^q\big( 
\bigoplus_{0}^{r} Sym^d(V_{k+1}^*)\big),
\end{equation}
$$ 1 \otimes \xi \in Det \otimes \big(\bigoplus_{0}^r 
Sym^d(V_{k+1}^*)\big).$$
The weight of $1\otimes 1 \in Det\otimes \bold{1}$ is
$\sum_{0}^{k} w_i$.
Since the polynomials $\{f_l\}$ do
not simultaneously vanish at $[1,0, \ldots, 0]\in \proj^k$,
one of the coefficients of $v_0^{*d}$ among the
polynomials $\{ f_l\}$ must be non-zero.
Similarly non-zero coefficients of $v_1^{*d}, \ldots, v_{k}^{*d}$
can be found among the polynomials $\{f_l\}$.
Therefore, the terms
\begin{equation}
\label{jayy}
 1\otimes (v_j^{*d} \otimes \cdots \otimes v_j^{*d})
\end{equation}
occur in (\ref{termmm}) and have 
weight $-qd\cdot w_j + \sum_{0}^{k}
 w_i$.

There are now two cases. First assume $\sum_{0}^{k} w_i >0$,
then $1\otimes 1$ has positive weight.
If $\sum_{0}^{k} w_i \leq 0$, there must exist $j$ such that
$w_j <0$. Let $w_j$ be the negative
weight of greatest absolute value.
Hence, for all $i$, if $w_i<0$, then $-w_j +w_i \geq 0$.
Finally, since $q>k+1$,
$$-qd \cdot w_j + \sum_{0}^{k} w_i >0.$$ 
The term (\ref{jayy}) therefore has positive weight.
The Numerical Criterion
implies $\xi$ is a stable point for the
$\bold{GL}_{k+1}$-action. 
\end{pf}

As a consequence of Lemma (\ref{geintth}), 
$U(k,r,d)/\bold{GL}_{k+1}$ exists as a quasi-projective
variety. Standard arguments show that the
space
$$M_{\proj^k}(\proj^r,d) \eqq U(k,r,d)/\bold{GL}_{k+1}$$
has the desired functorial properties.
Note: the family of maps 
\begin{equation}
\label{fammy}
(\pi:\cal{P} \rarr S, \ \mu:
\cal{P}\rarr \proj^r)
\end{equation}
may not be a Zariski locally
trivial $\proj^k$-bundle over $S$. A Galois cover
construction is required to obtain the canonical algebraic
morphism
\begin{equation}
\label{cann}
S \rarr M_{\proj^k}(\proj^r, d)
\end{equation} induced by the family 
(\ref{fammy}). Alternatively, one can
define a map to $M_{\proj^k}(\proj^r,d)$ locally
in the \'etale topology on $S$. The morphism (\ref{cann})
is then obtained via {\em descente}.
Since $M_{\proj^1}(\proj^r,d)$ and $M_{0,0}(\proj^r,d)$
coarsely represent the same functor, these spaces are
canonically isomorphic.

Since $U(k,r,d)$ is nonsingular, contained in the
stable locus, and the $GL_{k+1}$-action has finite stabilizers,
Luna's Etale Slice Theorem can be applied to conclude $M_{\proj^k}
(\proj^r,d)$ has finite quotient singularities (see [L]).
Luna's Theorem requires a characteristic zero hypothesis.

Finally, since $U(k,r,d)$ is  equivariant and contained in a 
G.I.T. stable locus, the group action
\begin{equation}
\label{proppp}
GL_{k+1} \times U(k,r,d) \rarr U(k,r,d)
\end{equation}
is a proper action. This is established in [MFK], Corollary (2.5). 
The properness of the action (\ref{proppp}) is needed in
section (\ref{genarg}).

\section{\bf{The Homomorphism $\lambda:Ch(k,r,d) \rarr Ch(k,r,1)$}}
\label{homlam}
Let
$\nu: \proj^r \rarr \proj^r$ be a regular
map satisfying $\nu^{*}(\oh_{\proj^r}(1)) \eqq \oh_{\proj^r}(d)$.
The map $\nu$ induces a canonical morphism
$\tau_\nu: \grass(\proj^k, \proj^r)\rarr M_{\proj^k}(\proj^r,d)$
by the following considerations.
Let $\pi: \cal{P} \rarr \bold{G}(\proj^k, \proj^r)$ be the
tautological $\proj^k$-bundle over the Grassmannian. 
Since 
\begin{equation}
\label{famm1}
\cal{P} \subset \grass(\proj^k,\proj^r) \times
\proj^r,
\end{equation} there is
a canonical projection $\eta: \cal{P} \rarr \proj^r$.
Let $\mu: \cal{P} \rarr \proj^r$ be determined
by $\mu= \nu \circ \eta$. The family
\begin{equation}
\label{famm}
(\pi: \cal{P} \rarr \bold{G}(\proj^k, \proj^r),\  
\mu: \cal{P} \rarr \proj^r)
\end{equation}
is a family over $\grass(\proj^k, \proj^r)$ 
of degree $d$ maps from $\proj^k$ to 
$\proj^r$. Since $M_{\proj^k}(\proj^r,d)$ is a coarse moduli
space, the family (\ref{famm}) induces
a morphism from the base to moduli:
$$\tau_{\nu}:\bold{G}(\proj^k, \proj^r) \rarr
M_{\proj^k}(\proj^r,d).$$
Let  $\tau_{\nu}^*$ be the ring homomorphism induced by pull-back:
$$\tau^*_{\nu}: Ch(k,r,d) \rarr Ch(k,r,1).$$ 
Since $M_{\proj^k}(\proj^r,d)$ has finite quotient singularities,
the pull-back map $\tau^*_{\nu}$ is well defined (see [V]).
\begin{pr}
\label{homomor}
The homomorphism 
 $\tau_{\nu}^*$ does not depend upon $\nu$ and is
a graded ring isomorphism.
\end{pr}
\noindent 
Let $\lambda: Ch(k,r,d) \rarr Ch(k,r,1)$ be the
ring isomorphism $\tau_{\nu}^*$ for any regular map $\nu$.
Theorem (\ref{mainn}) is a
consequence of Proposition (\ref{homomor}).

The proof of Proposition (\ref{homomor}) will
be undertaken in several steps. 
First the independence result will be established
in Lemma (\ref{indy5}).
A surjectivity Lemma will be also be proven
in this section. The injectivity
of $\tau_{\nu}^*$ will
be proven in section (\ref{genarg}).

\begin{lm}
\label{indy5}
The homomorphism 
 $\tau_{\nu}^*$ does not depend upon $\nu$.
\end{lm}
\begin{pf}
Let $U(r,r,d) \subset \bigoplus_{0}^{r} H^0(\proj^r,
\oh_{\proj^r}(d))$ be the Zariski open locus of
basepoint free $(r+1)$-tuples of polynomials as defined
in section
(\ref{construct}). There is a tautological morphism
$$\nu_U: U(r,r,d) \times \proj^r \rarr \proj^r.$$
The tautological family (\ref{famm1}) over the 
Grassmannian pulls-back to a tautological family $\cal{P}_U$ over
$$\grass(\proj^k, \proj^r) \times U(r,r,d).$$ 
$\cal{P}_U$ is equipped with a canonical projection 
$$\eta_U: \cal{P}_U \rarr U(r,r,d) \times \proj^r.$$
Let $\mu_U= \nu_U \circ \eta_U$. The map $\mu_U$ defines
a family of degree $d$ maps from $\proj^k$ to $\proj^r$
over $\grass(\proj^k, \proj^r) \times U(r,r,d)$. There
is an induced map 
$$\tau_{U}:\grass(\proj^k, \proj^r) \times U(r,r,d)
\rarr M_{\proj^k}(\proj^r,d).$$
The morphism $\tau_{\nu}$ is induced
by the composition of the inclusion
$$i_{\nu}:\grass(\proj^k, \proj^r)
\rarr \grass(\proj^k, \proj^r) \times [\nu]
\subset  \grass(\proj^k, \proj^r) \times U(r,r,d)$$
with $\tau_{U}$.
Hence, $\tau_{\nu}^*= i_{\nu}^* \circ \tau_U^*$
Since $U(r,r,d)$ is an open set in affine space,
$i_{\nu}^*= i_{\nu'}$ for any two maps
$[\nu], [\nu '] \in U(r,r,d)$.
\end{pf}

If $k=r$, then $\grass(\proj^r, \proj^r)$ is a point
and $\tau_{\nu}^*$ is surjective. Assume $k<r$.
Let $1\leq j \leq r-k$. Let $H_j \subset \proj^r$
be a linear subspace of codimension $k+j$.
Define an algebraic subvariety
$C(H_j)\subset M_{\proj^k}(\proj^r,d)$ by the following
condition. $C(H_j)$ is the set of maps that meet
$H_j$. $C(H_j)$ is easily seen to be an irreducible
subvariety of codimension $j$ in $M_{\proj^k}(\proj^r,d)$.
There is a natural $\bold{GL}_{r+1}$-action on $M_{\proj^k}
(\proj^r,d)$ obtained from the symmetries of $\proj^r$.
Let $\xi \in \bold{GL}_{r+1}$. Certainly
$$\xi( \ C(H_j)\ )= C(\ \xi(H_j)\ ).$$
Since $\bold{GL}_{r+1}$ is a connected rational
group, the class $\sigma_j$ of $C(H_j)$ in the Chow ring
$Ch(k,r,d)$ is well-defined (independent of $H_j$).
\begin{lm}
\label{abcd} The pull-back of the class $\sigma_j$ for
$1\leq j \leq r-k$ is determined by:
$$\tau_{\nu}^*(\sigma_j)= d^{k+j}\cdot \sigma_j.$$
\end{lm}

\begin{pf}
Let $\nu: \proj^r \rarr \proj^r$ be a fixed
morphism satisfying $\nu^*(\oh_{\proj^r}(1))\eqq 
\oh_{\proj^r}(d)$. Let $H_j\subset \proj^r$ be
a general (with respect to $\nu$) linear space.
By Bertini's Theorem,
$\nu^{-1}(H_j)$ is a nonsingular complete
intersection of $k+j$ hypersurfaces of degree $d$
in $\proj^r$. The set theoretic inverse image
$\tau_{\nu}^{-1}( C(H_j))$ is the set of $k$-planes of
$\proj^r$ meeting $\nu^{-1}(H_j)$. A simple
tangent space argument shows that 
the scheme theoretic inverse image
$\tau_{\nu}^{-1}(C(H_j))$ is generically reduced. Hence, 
$$\tau_{\nu}^* (\sigma_j) = [\tau_{\nu}^{-1}(C(H_j))]\in Ch(k,r,1).$$
It remains to determine $[\tau_{\nu}^{-1}(C(H_j))]\in Ch(k,r,1).$

Recall $\pi:\cal{P}\rarr \grass(\proj^k,\proj^r)$
is the tautological $\proj^k$-bundle over the Grassmannian.
Let $L$ be the Chern class of the line bundle
$\eta^*(\oh_{\proj^r}(1))$ on  $\cal{P}$.
The following equations hold:
$$\pi_*(L^{k+j})= \sigma_j,$$
$$\pi_*( (d\cdot L)^{k+j})= [\tau_{\nu}^{-1}(C(H_j))].$$
These equations imply $\tau_{\nu}^*(\sigma_j)=d^{k+j}\cdot 
\sigma_j.$ 
\end{pf}

Consider the $d=1$ case, $\grass(\proj^k, \proj^r)\eqq
M_{\proj^k}(\proj^r,1)$. There is a tautological
bundle sequence on  $\grass(\proj^k, \proj^r)$:
$$0 \rarr S \rarr \com^{r+1} \rarr Q \rarr 0.$$
$Q$ is a bundle of rank $r-k$. For
$1 \leq j \leq r-k$, let $c_j(Q) \in Ch(k,r,1)$
be the $j^{th}$ Chern class of $Q$. It is well known that
$$c_j(Q)=\sigma_j.$$
Also, the classes $c_j(Q)\in Ch(k,r,1)$ generate
$Ch(k,r,1)$ as ring. These facts can be found, for example,
in [F]. Therefore, the following Lemma
is a consequence of Lemma (\ref{abcd}). 
\begin{lm}
The  homomorphism
$\tau_{\nu}^*: Ch(k,r,d) \rarr Ch(k,r,1)$ is surjective.
\end{lm}
\noindent
In fact, the subring of $Ch(k,r,d)$ generated by
$\sigma_1, \ldots, \sigma_{r-k}$ surjects onto
$Ch(r,k,1)$ via $\tau_{\nu}^*$.

\section{\bf{Generation of $Ch(1,r,d)$}}
\label{genn}
In order to complete the proof of Proposition (\ref{homomor}),
results on the generation of $Ch(k,r,d)$ are needed.
In this section, a special argument in the $k=1$ case is
developed. In sections (\ref{eqchg})-(\ref{genarg}), a general
generation argument using the theory of equivariant Chow
groups is established. The general argument also
covers the $k=1$ case. The special method for the
$k=1$ case involves a natural stratification of
$M_{\proj^1}(\proj^r,d)$. Unfortunately, this stratification
does not easily generalize when $k>1$.

Let $0\leq j\leq r-1$. Let $\sigma_0\in Ch(1,r,d)$ be the
unit (the fundamental class). Let $\sigma_{j\neq 0}$ be the class 
defined in section (\ref{homlam}).

\begin{pr}
\label{gen}
The elements $\sigma_i \cdot \sigma_j$ ($0\leq i \leq j \leq r-1$) 
span a $\Q$-basis of $Ch(1,r,d)$.
\end{pr}

The proof of Proposition (\ref{gen}) uses the
3-pointed moduli space of maps $M_{0,3}(\proj^r,d)$.
Let $1,2,\infty \in \proj^1$ be three marked points.
There is a natural isomorphism:
\begin{equation}
\label{udef}
M_{0,3}(\proj^r,d)\eqq \proj(U)=U(1,r,d)/\com^*
\subset \proj(\bigoplus_0^r H^0(\proj^1, 
\oh_{\proj^1}(d)))
\end{equation}
where $U(1,r,d)$ is the basepoint free locus (see section 
(\ref{construct})).
An element of $\proj(U)$ corresponds to a degree $d$ map
from $\proj^1$ to $\proj^r$ with the
three markings $1,2,\infty\in \proj^1$.
A map $[\mu]\in  M_{0,3}(\proj^r,d)$ corresponds to
a point in $\proj(U)$ by identifying the three markings of $[\mu]$
with the points $1,2,\infty\in \proj^1$. A tangent space
argument shows this identification is an isomorphism
of schemes (both are non-singular varieties).

The proof of Proposition (\ref{gen})
is a refinement of the methods that 
appear in [P]. 
For $0\leq j \leq r-1$, let $H_j \subset \proj^r$
be a linear space of codimension $1+j$. For 
$0 \leq a,b \leq r-1$, let $C(H_a,H_b)\subset M_{0,0}
(\proj^r,d)$  be the subvariety of maps meeting $H_a$ and $H_b$
(where $H_a$ and $H_b$ are in general position). 
A simple argument shows  the
equation $$[C(H_a,H_b)]=\sigma_a \cdot \sigma_b$$
holds in $Ch(1,r,d)$.
Note: intersection with the hyperplane $H_0$ imposes
no condition on the maps. In particular,
$C(H_0, H'_0)= M_{0,0}(\proj^r,d)$.

\begin{lm}
\label{maxx}
Let $0\leq a,b \leq r-1$.
Assume $(a,b) \neq (r-1,r-1)$.
Let $H_a, H_b \subset \proj^r$ be linear spaces of
codimension $1+a, 1+b$
in general position.  
Let $H_{a+1}\subset H_a$, $H_{b+1}\subset H_b$ be 
linear spaces of codimension $1$. 
The natural map 
\begin{equation}
\label{mlem}
C(H_{a+1}, H_b)\ 
\cup \ C(H_a, H_{b+1})\ \cup \
C(H_0,H_{a}\cap H_{b})
\rarr 
C(H_a,H_b)
\end{equation}
yields a surjection on Chow groups
of proper codimension in $C(H_a, H_b)$.
If the linear spaces  $H_{a+1}$, $H_{b+1}$, or $H_a \cap H_b$ 
are empty, the corresponding cycle on the left in (\ref{mlem})
is taken to 
be empty. By the assumption $(a,b)\neq (r-1,r-1)$,
not all cycles are empty.
\end{lm}

\begin{pf}
Let $F$ be a hyperplane in general position
with respect to $H_a$ and $H_b$.
Let $\barr{N}=\barr{M}_{0,3}(\proj^r,d)$ and
$\barr{M}=\barr{M}_{0,0}(\proj^r,d)$. 
Let $N$, $M$ be the unbarred moduli spaces.
Let
$e_i:\barr{N} \rarr \proj^r$ be the natural
evaluation maps for the markings $1\leq i \leq 3$.
Let $$X= e_1^{-1}(F) \cap e_2^{-1}(H_a) \cap
e_3^{-1}(H_b).$$ 
$X$ is closed subvariety of $\barr{N}$.
The natural forgetful morphism
$\rho: X\rarr \barr{M}$ is proper. Also $\rho(X) \cap M=C(H_a,H_b)$.
Let $Z\subset C(H_a,H_b)$ be the open set of $\rho(X)$
corresponding to Kontsevich
stable maps satisfying the following conditions:
\begin{enumerate}
\item[(i.)] The domain curve is $\proj^1$.
\item[(ii.)] The map meets $H_a$ and $H_b$.
\item[(iii.)] The map does not pass through $F\cap H_a$,
     $F\cap H_b$, or $H_a \cap H_b$.
\end{enumerate}
Let $[\mu] \in Z$ be a element.
By condition (iii), the image $Im(\mu)\subset \proj^r$ can not be
contained in $F$, $H_a$, or $H_b$.
Moreover, by (iii), $\rho^{-1}(Z) \subset N$.
Hence, the map $\rho^{-1}(Z)\rarr Z$ 
has finite fibers. Since   $\rho^{-1}(Z)\rarr Z$ is
a proper morphism with finite fibers, it is a finite morphism.
Therefore, if $A_i(\rho^{-1}(Z)) \otimes \Q=0$, then
$A_i(Z) \otimes \Q=0$.

The set $\rho^{-1}(Z)\subset N \eqq \proj(U)$ (see (\ref{udef}) above)
is isomorphic to a quasi-projective variety 
in $\proj(\bigoplus_0^r H^0(\proj^1, 
\oh_{\proj^1}(d)))$.
The quasi-projective subvariety  $$\rho^{-1}(Z)\subset 
\proj(\bigoplus_0^r H^0(\proj^1, 
\oh_{\proj^1}(d)))$$ can be identified as follows.
Let $L_1\subset \proj(U)$ correspond
to maps $\mu:\proj^1 \rarr \proj^r$ satisfying $\mu(1)\in F$.
Let $L_2$, $L_{\infty}$ be the linear
subspaces in $\proj(U)$ where $\mu(2)\in H_a$, $\mu(\infty)\in H_b$.
Let $L_1\cap L_2 \cap L_{\infty}= I\subset \proj(U)$.
Let $D\subset I$ be the union of the 
three hypersurfaces of maps meeting
the linear spaces $F\cap H_a$, $F\cap H_b$, and 
$H_a\cap H_b$ respectively.
Since $(a,b) \neq (r-1,r-1)$,
  $F\cap H_a$ or $F\cap H_b$ is non-empty. Therefore,
$D\subset I$ is a subvariety of codimension $1$.
Then 
$$\rho^{-1}(Z)=  I \setminus D.$$
$I$ is an open set of a linear subspace of
 $\proj(\bigoplus_0^r H^0(\proj^1, 
\oh_{\proj^1}(d)))$. Since $D$ is of codimension $1$ in $I$, 
all the Chow groups of $\rho^{-1}(Z)$ of proper codimension vanish.
Hence all the Chow groups (tensor $\Q$) of $Z$ of proper
codimension also
vanish. 

By definition, $Z\subset C(H_a, H_b)$. 
The complement of $Z$ in $C(H_a, H_b)$ is the set of maps 
meeting $F\cap H_a$, $F\cap H_b$, or $H_a\cap H_b$.
Therefore, the complement of $Z$ in $C(H_a,H_b)$ is the union
of three cycles: 
\begin{equation}
\label{onion3}
C(F\cap H_a, H_b)\ \cup\ C(H_a,F\cap H_b)\ \cup\
C(H_0,H_a\cap H_b).
\end{equation}
Since the Chow groups  of $Z$ vanish in proper
codimension, the Chow groups of the union (\ref{onion3})
surject onto the Chow groups of $C(H_a, H_b)$ (in proper
codimension).
\end{pf}

\noindent A vanishing result is also required.
\begin{lm}
\label{four4}
Chow groups in proper codimension of 
$C(H_{r-1}, H'_{r-1})$ vanish.
\end{lm}

\begin{pf}
Let $F$ be a hyperplane in general position
with respect to two distinct points  $p=H_{r-1}$ and 
$q=H'_{r-1}$.
The notation $N\subset\barr{N}$, $M\subset \barr{M}$
of Lemma (\ref{maxx}) will be used.
Let $$X=e_1^{-1}(F) \cap e_2^{-1}(p)\cap e_3^{-1}(q).$$
Let $\rho: X \rarr \barr{M}$ be the proper forgetful morphism.
Again, $\rho(X) \cap M=C(p,q)$.
Let $Z\subset C(p,q)$ be the open set of $\rho(X)$
corresponding to Kontsevich
stable maps satisfying the following conditions:
\begin{enumerate}
\item[(i.)] The domain curve is $\proj^1$.
\item[(ii.)] The map meets $p$ and $q$.
\end{enumerate}
Note $F\cap p$, $F\cap q$, and $p\cap q$ are empty.
By these conditions on $Z$, the map $\rho^{-1}(Z)\rarr Z$ is
finite and proper. 
Therefore, if $A_i(\rho^{-1}(Z)) \otimes \Q=0$, then
$A_i(Z) \otimes \Q=0$. Also, $\rho^{-1}(Z) \subset N$.

The quasi-projective subvariety  $$\rho^{-1}(Z)\subset 
\proj(\bigoplus_0^r H^0(\proj^1, 
\oh_{\proj^1}(d)))$$ can be identified as follows.
Let $L_1\subset \proj(U)$ correspond
to maps $\mu:\proj^1 \rarr \proj^r$ satisfying $\mu(1)\in F$.
Let $L_2$, $L_{\infty}$ be the linear
subspaces in $\proj(U)$ where $\mu(2)\in p$, $\mu(\infty)\in q$.
Let $L_1\cap L_2 \cap L_{\infty}= I\subset \proj(U)$.
Then 
$$\rho^{-1}(Z)=  I $$
$I$ is an open set of a linear subspace of
 $\proj(\bigoplus_0^r H^0(\proj^1, 
\oh_{\proj^1}(d)))$. 
Let $\barr{I}$ be the closure of $I$. It will be shown
that $\barr{I}\setminus I$ has codimension $1$ in $\barr{I}$.
The Chow groups of $\rho^{-1}(Z)$ of proper codimension therefore
vanish.
Hence all the Chow groups of $Z$ of proper
codimension also
vanish. 

Let $[A_0, \ldots, A_r]$ be homogeneous coordinates on $\proj^r$.
Let  $$F=(A_0-A_r),\  p=[1,0, \ldots,0], \  q=[0,\ldots,0,1].$$
Let $[S,T]$ be homogeneous coordinates on
$\proj^1$. Let $1,2,\infty\in \proj^1$ be the points
$[1,1]$, $[1,0]$, $[0,1]$ respectively. 
An element $[\mu]\in \proj(\bigoplus_0^r H^0(\proj^1, 
\oh_{\proj^1}(d)))$ is given by
an $r$-tuple of degree $d$ homogeneous 
polynomials in $S$ and $T$ :  $[f_0, \ldots, f_r]$.
The element $[\mu]$ is in $I$ if and only if
\begin{enumerate}
\item[(i.)] $f_0,\ldots, f_r$ span a basepoint free
            linear system on $\proj^1$.
\item[(ii.)] $f_0(1,1)=f_r(1,1)$.
\item[(iii.)]  $T$ divides $f_1, \ldots, f_r$.
\item[(iv.)] $S$ divides $f_0, \ldots, f_{r-1}$.
\end{enumerate}
The additional condition
$$S \ \  divides\ \  f_r$$ is a codimension
$1$ condition contained in the set $\barr{I}\setminus I$.
Hence $\barr{I}\setminus I$ has codimension $1$ in $I$.
\end{pf}
\noindent
Repeated application of Lemma (\ref{maxx}) shows the ring
$Ch(1,r,d)$ is generated (as a $\Q$-vector space) by the
classes $[C(H_a, H_b)]$ and the Chow groups of
$C(H_{r-1}, H'_{r-1})$. Lemma (\ref{four4}) shows
the Chow groups of $C(H_{r-1}, H'_{r-1})$ vanish in
proper codimension. Hence the classes
$[C(H_a,H_b)]=\sigma_a \cdot \sigma_b$ generate $Ch(1,r,d)$.

Via the classical
Schubert calculus,
the classes $\sigma_a \cdot \sigma_b$ for $0\leq a,b \leq r-1$
are easily seen to span a {\em basis} of the Chow ring of the
linear
Grassmannian $Ch(1,r,1)$. Consider the ring homomorphism
$$\tau^*_{\nu}: Ch(1,r,d) \rarr Ch(1,r,1)$$
defined in section (\ref{homlam}).
By Lemma (\ref{abcd}), 
$$\tau^*_{\nu}(\sigma_0)=\sigma_0,$$
$$\forall a> 0, \ \ \tau^*_{\nu}(\sigma_a)= d^{1+a}\sigma_a,$$
$$\forall a,b> 0, \ \ \tau^*_{\nu}(\sigma_a
\cdot \sigma_b)= d^{2+a+b}\sigma_a\cdot \sigma_b,$$
Therefore, the elements $\sigma_a \cdot \sigma_b$ for 
$0\leq a,b \leq r-1$ are independent in $Ch(1,r,d)$.
Since generation was established above, Proposition (\ref{gen})
is proven.
In case $k=1$, the injectivity of $\tau_{\nu}^*$ has been
proven.

\section{Equivariant Chow Groups}
\label{eqchg}
Let $G$ be a group. 
Let $G\times X \rarr X$ be a left group action.
In topology, the $G$-equivariant cohomology of 
$X$ is defined as follows. Let $EG$ be a contractable
topological space equipped with a free left $G$-action and 
quotient $EG/G=BG$.
Consider the left action of $G$ on $X\times EG$ defined
by:
$$g(x,b)= (g(x), g(b)).$$
$G$ acts freely on $X\times EG$. Let
$X\times_{G} EG$ be the (topological) quotient. 
The $G$-equivariant cohomology of 
of $X$, $H_G^*(X)$, is defined by:
$$H_G^*(X) = H^*_{sing}(X\times_{G} EG).$$
If $X$ is a a locally trivial principal $G$-bundle,
then $X\times_{G} EG$ is a locally trivial
fibration of $EG$ over the quotient $X/G$.
In this case, $X\times_{G} EG$ is homotopy equivalent
to $X/G$ and
$$H_G^*(X) = H^*_{sing}(X\times_{G} EG) \eqq H^*_{sing}(X/G).$$
For principal bundles, computing the
equivariant cohomology ring is equivalent
to computing the cohomology of the quotient.

There is an analogous equivariant theory of Chow groups
developed by D. Edidin, W. Graham, and B. Totaro in
[EG], [T]. Let $G$ be a linear algebraic group.
Let $G\times X \rarr X$ be a linearized, algebraic $G$-action.
The algebraic analogue of $EG$ is attained by
approximation. Let $V$ be a $\com$-vector space. Let
$G\times V \rarr V$ be an
algebraic representation of of $G$. 
Let $W\subset V$ be a $G$-invariant open set satisfying:
\begin{enumerate}
\item[(i.)] The complement of $W$ in $V$ is of codimension greater than q.
\item[(ii.)] $G$ acts on $W$ with trivial stabilizers.
\item[(iii.)] There exists a geometric quotient $W\rarr W/G$.
\end{enumerate}
$W$ is an approximation of $EG$ up to codimension $q$.
By (iii) and the assumption of
linearization, a geometric quotient $X\times _{G} W$ exits as
an algebraic variety. Let $d=dim(X)$, $e=dim( X\times _{G} W)$.
The equivariant Chow groups are defined by:
\begin{equation}
\label{defff} 
A^{G}_{d-j}(X)= A_{e-j}(X\times _{G} W)
\end{equation}
for $0\leq j \leq q.$
An argument is required to check these equivariant
Chow groups are well-defined (see [EG]).
The basic functorial properties of equivariant
Chow groups are  
established in [EG]. In particular, if $X$ is
nonsingular, there is a natural intersection
ring structure on $A_i^{G}(X)$.

Let $Z$ be a variety of dimension $z$.
For notational convenience, a superscript will
denote the Chow group codimension:
 $$A^{G}_{z-j}(Z) = A^j_G(Z), \ A_{z-j}(Z)=A^j(Z).$$
In particular, equation (\ref{defff}) becomes:
$$\forall\  0\leq j \leq q, \ \
A_{G}^{j}(X)= A^j(X\times _{G} W).$$
The following result of [EG] 
will be used in section (\ref{genarg}).
\begin{pr}
\label{dane}
Let $\com$ be the ground field of complex numbers.
Let $X$ be a quasi-projective variety. Let $G$ be a reductive
group.
Let $G\times X \rarr X$ be a linearized, proper, $G$-action
with finite stabilizers. Let $X\rarr X/G$ be a
quasi-projective
geometric quotient. Then, there are natural 
isomorphisms for all $j$:
$$ A^{j}_G(X) \otimes \Q \eqq A^{j}(X/G) \otimes \Q.$$
\end{pr}

\section{\bf{The Chow Ring of the Grassmannian 
and} $A_*^{GL}(pt)$}
\label{grdwk}
Let $0 \rarr S \rarr \com^{r+1} \rarr Q \rarr 0$ be
the tautological sequence on $\grass(\proj^k, \proj^r)$.
The following presentation of the Chow ring will be used
in section (\ref{genarg}). Let
$$c_1, \ldots, c_{k+1} \in Ch(k,r,1)$$ be the Chern
classes of the rank $k+1$ bundle $S$.
These classes generate $Ch(k,r,1)$.
Let $$c(S)= 1+ c_1 \ t+ c_2 \ t^2 + \cdots + c_{k+1} 
\ t^{k+1},$$
$${1 \over c(S)} = 1+p_1(c_1) \ t+ p_2(c_1,c_2) 
\ t^2+ p_3(c_1,c_2,c_3)\   t^3+
\cdots $$ where the latter is 
the formal inverse in the formal power series ring
$\com[c_1, \ldots, c_{k+1}][[t]]$.
The ideal of relations among $c_1, \ldots, c_{k+1}$ in $Ch(k,r,1)$
is generated by
the polynomials $$\{p_j \ | \ j> r-k\}.$$ Geometrically,
these relations are obtained from the vanishing of the $j^{th}$
Chern class  of the rank $r-k$ bundle $Q$ for $j>r-k$.

In section (\ref{genarg}), a basic result on push-forwards is needed.
\begin{lm}
\label{bozo}
Let $\pi:\proj(S) \rarr \grass(\proj^k,\proj^r)$
be the canonical projection. Let $\oh_{\proj(S)}(1)$
be the canonical line bundle on $\proj(S)$. Then, for $l\geq k$,
\begin{equation}
\label{ecc}
\pi_*( \ c_1^l(\oh_{\proj(S)}(1))\ )= p_{l-k} \in Ch(k,r,1).
\end{equation}
\end{lm}
\begin{pf}
Let $\xi= c_1(\oh_{\proj(S)}(1))$.
Certainly, $\pi_*(\xi^k)= 1= p_0$.
The equation
$$ \xi^{k+1} + c_1\  \xi^k+ \cdots +c_k \ \xi+  c_{k+1}=0$$
recursively yields (\ref{ecc}).
\end{pf}

The equivariant Chow ring $A_*^{\bold{GL}_{k+1}}(pt)$ is computed
to motivate the construction in (\ref{genarg}).
The notation of section (\ref{construct}) will be used here.
Let $V_{k+1}$ be a fixed $k+1$- dimensional complex
vector space such that $\proj(V_{k+1})=\proj^k$,
Let $U(k,n,1) \subset \bigoplus_{0}^{n}V^*_{k+1}$
be the basepoint free locus.
The codimension of the complement of
$U(k,n,1)$ is easily found to be $n-k+1$. 
$\bold{GL}(V_{k+1})$ acts on $U(k,n,1)$ with trivial stabilizers.
As determined in section (1), there is a geometric quotient
$$U(k,n,1)/\bold{GL}(V_{k+1}) \eqq \grass(\proj^k,\proj^n).$$
By the definition of the equivariant Chow ring,
$$A^{j}_{\bold{GL}_{k+1}} (pt)= A^{j}(\grass(\proj^k,\proj^n))$$
for $0\leq j \leq n-k$.
By the presentation of the Chow ring of  $\grass(\proj^k,\proj^n)$
given above, the relations
among the generators $c_1, \ldots, c_{k+1}$
start in codimension $n-k+1$. Hence, $A^{*}_{\bold{GL}_{k+1}}(pt)$
is freely generated (as a ring) by $c_1, \ldots, c_{k+1}$ where
$c_j \in A^{j}_{\bold{GL}_{k+1}}(pt)$.

\section{\bf{The Generation Argument}}
\label{genarg}
Again, let $U(k,n,1) \subset \bigoplus_{0}^{n}V^*_{k+1}$
be the basepoint free open set. As $n\rarr \infty$, 
$U(k,n,1)$ approximates
$E\bold{GL}_{k+1}$. 
By the definitions,
$$A^{j}_{\bold{GL}_{k+1}}(U(k,r,d)) \eqq
 A^{j}\big( U(k,r,d) \times_{\bold{GL}_{k+1}} U(k,n,1)\big)$$
for $0 \leq j \leq n-k$.
Recall
$$U(k,r,d) \subset \bigoplus_{0}^{r} Sym^d(V_{k+1}^*)$$
is the basepoint free locus.
There is a natural $\bold{GL}(V_{k+1})$-equivariant
 open inclusion,
$$U(k,r,d) \times U(k,n,1) \subset  \bigoplus_{0}^{r} Sym^d(V_{k+1}^*)
 \times U(k,n,1),$$
which yields an open inclusion
$$U(k,r,d) \times_{\bold{GL}_{k+1}}
 U(k,n,1) \subset  \bigoplus_{0}^{r} Sym^d(V_{k+1}^*)
 \times_{\bold{GL}_{k+1}} U(k,n,1).$$
Let $0 \rarr S \rarr \com^{n+1} \rarr Q \rarr 0$
be the tautological sequence on $\grass(\proj^k,\proj^n)$.
It is routine to verify
$$\bigoplus_{0}^{r} Sym^d(V_{k+1}^*)
 \times_{\bold{GL}_{k+1}} U(k,n,1) = \bigoplus_{0}^{r} Sym^d(S^*)$$
where the latter is the total space of the bundle
$\bigoplus_{0}^{r} Sym^d(S^*)$ over $\grass(\proj^k,\proj^n)$.
Let $D$ be the complement of $U(k,r,d) \times_{\bold{GL}_{k+1}}
U(k,n,1)$ in $\bigoplus_{0}^{r}Sym^d(S^*)$.

The Chow ring of $\bigoplus_{0}^{r} Sym^d(S^*)$ is isomorphic to
$Ch(k,n,1)$ via pull back. 
Let $dim$ be the dimension of the variety
$ \bigoplus_{0}^{r} Sym^d(S^*)$.
Let 
$$i_D: A_{dim-j}(D) \rarr A^j(\bigoplus_{0}^{r}
Sym^d(S^*))$$
be the map obtained by the inclusion 
$D\subset \bigoplus_{0}^{r} Sym^d(S^*)$.
There are exact sequences
of Chow groups
$$ A_{dim-j}(D) \rarr A^{j}(\bigoplus_{0}^{r}Sym^d(S^*)) \rarr
A^{j}\big(U(k,r,d) \times_{\bold{GL}_{k+1}}
U(k,n,1)\big)\rarr 0.$$
Let $c_1, \ldots, c_{k+1}$ be the classes of
$Ch(k,n,1)$ defined in section (\ref{grdwk}).
\begin{lm}
\label{bilt}
$p_j(c_1, \ldots, c_{k+1})
\in Im(i_D) \subset A^{j}(\bigoplus_{0}^{r}Sym^d(S^*))$
for all $j>r-k$.
\end{lm}
\begin{pf}
The proof involves an auxiliary construction.
Let $\pi:\proj(S) \rarr \grass(\proj^k,\proj^n)$
be the projective bundle associated to $S$.
Let $T= \pi^{*}\big( \bigoplus_{0}^{r}
 Sym^d(S^*) \big)$. Denote the total
space of the bundle $T$ also by $T$.
There is a commutative diagram.
\begin{equation*}
\begin{CD}
 T@>>>  \proj(S) \\
@V{\pi}VV @V{\pi}VV \\
\bigoplus_{0}^{r}Sym^d(S^*) @>>> \grass(\proj^k,\proj^n) \\
\end{CD}
\end{equation*}
There is a tautological rational evaluation map
$$\gamma: T \ - \ - \ \rarr \proj^r.$$
A point $\tau \in T$ is a triple 
$$\tau=(v, V\subset \com^{n+1}, (f_0, f_1, \ldots, f_r))$$
where $v$ is element of the $k$-dimensional projective
space $\proj(V)$ and $f_i \in Sym^d(V^*)$.
The rational map $\gamma$ is obtained by
$$\gamma (\tau)= [f_0(v), \ldots, f_r(v)].$$
Let $\barr{D}$ be the set of elements $\tau\in T$
such that all the $f_i$ vanish at $v$. $\barr{D}$ is
the  locus where $\gamma$ is undefined. The important
fact is
$$\pi(\barr{D})= D \subset \bigoplus_{0}^{r}Sym^d(S^*)$$
and $\pi: \barr{D} \rarr D$ is a projective, birational
morphism.
The Lemma will be proved by finding the
class of $\barr{D}$ in $A_*(T)$ and pushing forward.

Let $L$ be the line bundle $\oh_{\proj(S)}(d)$ on
$\proj(S)$. Let $L$ also denote the pull-back of
  $\oh_{\proj(S)}(d)$ to $T$. 
The rational map $\gamma$ is obtained from
$r+1$ tautological sections of $L$ on $T$.
There is a natural equivalence
$$ H^0(\grass(\proj^k,
\proj^n), Sym^d(S^*)) \eqq H^0(\proj(S),L).$$
Also, there is a natural inclusion 
$$ H^0(\grass(\proj^k,\proj^n),\ 
Sym^d(S^*) \otimes \bigoplus_{0}^{r} Sym^d(S)\ ) \subset H^0(T,L).$$
Since the bundle $Sym^d(S^*) \otimes Sym^d(S)$ has a
canonical identity section,
the bundle
$Sym^d(S^*)\otimes  \bigoplus _{0}^{r} Sym^d(S))$ has
$r+1$ canonical sections. 
It is straightforward to verify these $r+1$ sections
$z_0, \ldots, z_r$ of $L$ on $T$
yield the rational map $\gamma$. The  base locus
$\barr{D}$ is the common zero locus of the
sections $z_0, \ldots, z_{r}$. In fact, $\barr{D}$ is
a nonsingular variety of pure codimension $r+1$.
Explicit equations show $\barr{D}$ is 
nonsingular complete intersection. Hence
$[\barr{D}]= c_1(L)^{r+1} \in A_*(T)$.

Certainly $\pi_*(\ c_1(L)^{r+1}\ ) \in Im (i_D)$.
Also, for all $\alpha \geq 0$, 
$$\pi_*(\ c_1(L)^{r+1+\alpha}\ )= 
\pi_*([\barr{D}] \cap c_1(L)^{\alpha})\in Im(i_D).$$
It remains to compute $\pi_*(c_1(L)^{r+1+\alpha}) 
\in A_*(\bigoplus_{0}^{r} Sym^d(S^*))$.
But since push-forward commutes with flat pull-back, it suffices
to consider $c_1(L)^{r+1+\alpha} \in A_*(\proj(S))$ and compute
$\pi_*(c_1(L)^{r+1+\alpha})\in Ch(k,n,1)$. By Lemma (\ref{ecc}),
since $c_1(L)= d\cdot c_1(\oh_{\proj(S)}(1))$,
$$  \pi_*(c_1(L)^{r+1+\alpha})= d^{r+1+\alpha} \cdot
p_{r-k+1+\alpha}.$$
Hence $p_{j} \in Im(i_D)$ for all $j>r-k$.
\end{pf}

By Lemma (\ref{bilt}) and the presentation of $Ch(k,r,1)$
given in section (\ref{grdwk}), the following inequality is obtained:
\begin{equation}
\label{donn}
dim_{\Q} \ A^{j}_{\bold{GL}_{k+1}}(U(k,r,d)) \leq dim_{\Q}\
A^{j}(\grass(\proj^k,\proj^r)).
\end{equation}
Recall $\bold{GL}_{k+1} \times U(k,r,d) \rarr U(k,r,d)$
is a proper group action with finite stabilizers and
geometric quotient $M_{\proj^k}(\proj^r,d)$. Hence,
by Proposition (\ref{dane}) and inequality (\ref{donn}),
$$dim_{\Q} \ A^{j}(M_{\proj^k}(\proj^r,d) \leq dim_{\Q}
\ A^{j}(\grass(\proj^k,\proj^r)).$$
The surjectivity of $\tau_{\nu}^*:Ch(k,r,d) \rarr Ch(k,r,1)$
obtained in section (\ref{homlam}) implies
$$dim_{\Q} \ A^{j}(M_{\proj^k}(\proj^r,d) \geq dim_{\Q}
\ A^{j}(\grass(\proj^k,\proj^r)).$$
Therefore $\tau_{\nu}^*$ is injective.
The proofs of Proposition (\ref{homomor})
and Theorem (\ref{mainn}) are complete.
Since the subring of  $Ch(r,k,d)$ generated by the
classes $\sigma_1, \ldots, \sigma_{r-k}$ surjects onto
$Ch(r,k,1)$ via $\tau_{\nu}^*$, $Ch(r,k,d)$ is generated
(as a ring) by the classes $\sigma_1, \ldots, \sigma_{r-k}$.

\noindent
Department of Mathematics \\
University of Chicago\\
5743 S. University Ave\\
Chicago, IL 60637 \\
rahul@@math.uchicago.edu
\end{document}